\documentclass[11pt]{article}
\usepackage{amssymb,amsmath}
\usepackage[english]{babel}
\newtheorem{theo}{Theorem}[section]
\newtheorem{lema}[theo]{Lemma}
\newtheorem{prop}[theo]{Proposition}
\newtheorem{nota}[theo]{Remark}
\newtheorem{definicion}[theo]{Definition}
\newtheorem{corolary}[theo]{Corolary}

\newcommand{\bq}{\begin{equation}}
\newcommand{\eq}{\end{equation}}
\newcommand{\ba}{\begin{array}}
\newcommand{\ea}{\end{array}}

\newcommand{\racion}[2]{\mbox{\small$\frac{{#1}}{{#2}}$} }
\newcommand{\refe}[1]{(\ref{#1})}
\newcommand{\half}{\frac{1}{2}}
\newcommand{\tsig}{\tilde{\sigma}}
\newcommand{\ttau}{\tilde{\tau}}
\newcommand{\dst}{\displaystyle}

\newcommand\be{\begin{enumerate}}
\newcommand\ee{\end{enumerate}}
\newcommand\bi{\begin{itemize}}
\newcommand\ei{\end{itemize}}

\newcommand{\q}{\varsigma}
\newcommand{\pe}[2]{\langle#1,#2\rangle}

\newcommand\RR{\mathbb{R}}
\newcommand\CC{\mathbb{C}}

\newcommand\bd{\begin{definicion}{\bf }}
\newcommand\ed{\end{definicion}}
\newcommand\bl{\begin{lema}{\bf }}
\newcommand\el{\end{lema}}
\newcommand\bp{\begin{prop}{\bf }}
\newcommand\ep{\end{prop}}
\newcommand\bt{\begin{theo}{\bf }}
\newcommand\et{\end{theo}}
\newcommand\bdm{\begin{proof}}
\newcommand\edm{\end{proof}}
\newcommand\bn{\begin{nota}{\bf }}
\newcommand\en{\end{nota}}
\newcommand\bc{\begin{corolary}{\bf }}
\newcommand\ec{\end{corolary}}

\newenvironment{proof}{\noindent{\bf Proof:\ }}{\hfill
\fbox \par\vspace{2ex}}

\newfont{\got}{eufm10 scaled \magstep1}
\newcommand{\g}[1]{\mbox{\got #1}}

\title{Factorization of the hypergeometric-type difference
equation \\ on the non-uniform lattices:  dynamical algebra}
\date{\today}
\author{R. \'Alvarez-Nodarse${}^\dag{}^\ddag$,
N. M. Atakishiyev${}^\S$ and R. S.
Costas-Santos${}^*$\\[5mm]
\small ${}^\dag$ Departamento de An\'alisis Matem\'atico.\\
\small Universidad de Sevilla. Apdo. 1160, E-41080 Sevilla, Spain\\
\small${}^\ddag$ Instituto Carlos I de F\'{\i}sica Te\'orica y
Computacional, \\
\small Universidad de Granada, E-18071 Granada, Spain\\
\small${}^\S$ Instituto de Matem\'aticas, UNAM,  Apartado Postal 273-3, \\
\small C.P. 62210 Cuernavaca, Morelos, M\'exico\\
\small${}^*$ Departamento de Matem\'aticas,
E.P.S., Universidad Carlos III de Madrid.\\
\small Ave. Universidad 30, E-28911, Legan\'es, Madrid, Spain \\
 }

\begin{document}
\maketitle
\begin{abstract}
\noindent
We argue that one can factorize the difference equation of
hypergeometric type on the non-uniform lattices in general
case. It is shown that in the most cases of $q$-linear spectrum
of the eigenvalues this directly leads to the dynamical symmetry
algebra $su_q(1,1)$, whose generators are explicitly constructed
in terms of the difference operators, obtained in the process of
factorization.  Thus all models with the $q$-linear spectrum (some of them,
but not all, previously considered in a number of publications)
can be treated in a unified form.
\end{abstract}

\section{Introduction and preliminaries}

In this paper we continue the study, started in \cite{alco1}, on the
factorization of the hyper\-geo\-me\-tric-type  difference equation on
the non-uniform lattices, i.e., of the equation \cite{niuv1}
\begin{equation}
\begin{array}{c}
\dst \sigma(s) \frac{\Delta}{\Delta x(s-\frac{1}{2})}
\left[\frac{\nabla y(s)}{\nabla x(s)}\right] + \tau(s)
\frac{\Delta y(s)}{\Delta x(s)} + \lambda y(s) =0, \\[5mm]
\sigma(s)=\tilde{\sigma}(x(s)) - \frac{1}{2}\tilde{\tau}(x(s))
\Delta x\left(\mbox{$s-\frac{1}{2}$}\right), \quad
\tau(s)=\tilde{\tau}(x(s)),
\end{array}
\label{eqdif}
\end{equation}
where $\Delta y(s):= y(s+1)-y(s)$, $\nabla y(s):= y(s)-y(s-1)$,
$\tsig(x(s))$ and $\ttau(x(s))$ are polynomials in $x(s)$ of degree
at most 2 and 1, respectively, and $\lambda$ is a constant (see also
\cite{nisuuv}).
The difference equation \refe{eqdif} has polynomial solutions
$P_n(x(s);q)\!:=$  $ P_n(s;q)$ of the hypergeometric type if and
only if the lattice $x(s)$ has the form \cite{niuv2,atrasu}
\begin{equation}
x(s)= c_1(q) q^{s}+c_2(q)q^{-s}+c_3(q)=
c_1(q)[q^{s}+q^{-s-\mu}]+c_3(q), \label{red-gen}
\end{equation}
where $c_1$, $c_2$, $c_3$ and $q^{\mu}:=c_1/c_2$ are constants which,
in general, depend on $q$. An important special case of the lattice
$x(s)$ is the $q$-linear lattice, which is obtained from \refe{red-gen}
by assuming that either $c_1(q)$ or $c_2(q)$ vanishes.

The polynomial solutions of the difference equation \refe{eqdif}
correspond to the following expression \cite{nisuuv} for its eigenvalues
$\lambda_n(q)$:
\bq\label{lam1} \ba{c}
\lambda_n(q) = C_1q^n+C_2q^{-n}+C_3,\\[5mm]
\dst C_1 = \frac{1}{2(1-q)} \left(\widetilde{\tau}'
+ \frac{\widetilde{\sigma}''}{k_q}\right),
\quad C_2 = \frac{1}{2(1-q^{-1})}\left(\widetilde{\tau}' -
\frac{\widetilde{\sigma}''}{k_q}\right), \quad
C_3 = -\frac{\widetilde{\sigma}''(1+q)}{2k_q(1-q)} -
\frac{\widetilde{\tau}'}{2},
\ea
\eq
where $\widetilde{\tau}'$ and $\widetilde{\sigma}''/2$ are the
coefficients of $x(s)$ and $x^2(s)$ in the Taylor expansion for
$\ttau(x(s))$ and $\tsig(x(s))$, respectively, i.e., $\ttau(x(s))=
\widetilde{\tau}' x(s)+ \widetilde{\tau}(0)$, and $\tsig(x(s))=
\widetilde{\sigma}''/2\,x^2(s)+ \widetilde{\sigma}'(0)x(s)+
\widetilde{\sigma}(0)$.

Observe that the coefficients $C_1$ and $C_2$ of the $q^n$ and
$q^{-n}$ terms, respectively, are fixed by the functions $\sigma$
and $\tau$ in \refe{eqdif}, and so is the product $C_1C_2$. In what
follows we denote by $L_q$ the value of $C_1C_2 =
\left(\left({\widetilde{\sigma}''}/{k_q}\right)^2-(\widetilde{\tau}')^2
\right)/{4k_q^2}$.\\
The sequence $\left\{\lambda_n(q) \right\}$ satisfies the following
three-term recurrence relation (TTRR)
\bq\label{ttrr-lam}
\lambda_{n+2}(q)-(q+q^{-1})\,\lambda_{n+1}(q)+ \lambda_n(q)=
\frac12({\tilde\tau' k_q^2-\tilde\sigma''[2]_q})=C.
\eq
Conversely, if $\left\{\lambda_n(q)\right\}$ satisfy the TTRR
\refe{ttrr-lam}, then it has the form $\lambda_n(q) = C_1'q^n +
C_2'q^{-n}+ C_3'$. Obviously, having used the initial conditions
$\lambda_0(q) = 0$ and $\lambda_1(q) = -\widetilde\tau'$, one recovers the
expression \refe{lam1}.\\
It is well known \cite{nisuuv} that under certain conditions the
polynomial solutions of \refe{eqdif} are orthogonal. For example, if
$\sigma(s) \rho(s) x^{k}(s-\frac{1}{2}) \Big|_{s=a,b} = 0$, for all
$k=0,1,2,\dots$, then the polynomial solutions possess a discrete
orthogonality property \bq\label{dis-ort} \sum_{s=a}^{b-1}P_n(s;q)\,
P_m(s;q)\, \rho(s) \nabla x_1(s)= d_n^2(q) \delta_{n,m}, \eq
where the 
weight function $\rho(s)$ is a solution
of the Pearson-type difference equation \cite{nisuuv}
\bq\label{pearson}
\frac{\Delta}{\Delta x(s-\frac{1}{2})}
\left[\sigma(s) \rho(s)\right]= \tau(s) \rho(s)\quad\mbox{or}
\quad \sigma(s+1)\rho(s+1)=\sigma(-s-\mu)\rho(s).
\eq
If the lattice $x(s)$ is a $q$-linear lattice, i.e., $x(s)= c\,q^{\pm s}
+ c_3$, then the $\sigma(-s-\mu)$ in \refe{pearson} should be substituted
by $\sigma(s)+\tau(s)\Delta x(s-1/2)$. A more detailed information on
orthogonal polynomials on the non-uniform lattices can be found in
\cite{nisuuv,atrasu,gara,kost,alme}.
\\
In \cite{alco1} it has been shown that one can factorize the
Nikiforov-Uvarov equation \refe{eqdif} with the aid of raising and
lowering operators, which can be constructed for solutions of this
equation. In this paper we wish to make one step further by studying
the dynamical symmetry algebra for the hypergeometric-type
difference equation \refe{eqdif} on the non-uniform lattices
\refe{red-gen}. Our approach is essentially based on the simple
observation, formulated in \cite{ata1}: In order to factorize an
arbitrary difference equation, one should express it explicitly in
terms of the shift (or displacement) operators $\exp
(a\,\frac{d}{ds})$, which are defined as $\exp
(a\,\frac{d}{ds})\,f(s) = f(s+a)$, $a$ is some constant. For
example, in the case of the equation \refe{eqdif} this corresponds
to the substitutions $\Delta = \exp (\frac{d}{ds}) - 1$ and $\nabla
= 1 - \exp (-\frac{d}{ds})$. This procedure converts a difference
equation into an eigenvalue problem for a difference operator,
represented by a linear combination of some shift operators (with
coefficients, which depend polynomially on the variable $s$). Since
each term of this linear combination is readily factorizable
(because $\exp{(\alpha + \beta)\,A} =
\exp{\alpha\,A}\,\exp{\beta\,A}$ for an arbitrary operator $A$), the
factorization of the whole linear combination, which represents the
initial difference equation, becomes
straightforward.\\
Inspired by the appearance of Macfarlane's \cite{mac} and
Biedenharn's \cite{bie} important constructions of
$q$-a\-na\-lo\-gues of quantum harmonic oscillator, this technique
of factorization of difference equations was later employed in a
number of publications \cite{atsu2}--\cite{atfrwo} in order to
study group theoretic properties of the various well-known families
of orthogonal polynomials, which can be viewed as $q$-extensions of
the classical Hermite polynomials. So our purpose here is to
formulate a unified approach to deriving
all of these results, which correspond to the $q$-linear spectrum. \\
An important aspect to observe at this point is that we shall mainly
(except for the examples in subsection 4.2) confine our attention to
those families of $q$-polynomials, which satisfy discrete
orthogonality relation of the type \refe{dis-ort}. The explanation
of such preference is that the factorization of difference equations
for instances of $q$-polynomials with continuous orthogonality
property has been already thoroughly studied in
\cite{atsu2}--\cite{atfrwo}. Observe also that our approach still
remains valid in the limit as $q\to 1$; so classical counterparts of
$q$-polynomials, which will be discussed in this paper, are in fact
incorporated as appropriate limit cases. But the reader who desires
to know more about the factorization in the cases of classical
orthogonal polynomials (such as the Kravchuk, Charlier, Meixner,
Meixner--Pollaczek, and Hahn) may be referred to \cite{atjanawo,ban2}
and references therein. \\
The paper is organized as follows. In section \ref{fac-op} we
associate with each family of $q$-polynomials a ``$q$-Hamiltonian''
$\g H(s;q)$ (via the second-order difference equation) and construct
two difference operators $a(s;q)$ and $b(s;q)$, which factorize the
operator $\g H(s;q)$. Our main results are given in section
\ref{alg-din}: they are formulated in Theorems \ref{el-teo-q-lin} --
which gives a {\em simple} necessary and sufficient condition that
the $q$-Hamiltonian $\g H(s;q)$ admits the factorization in terms of
the operators $a(s;q)$ and $b(s;q)$, which satisfy the relation
$a(s;q)b(s;q)-q^\gamma b(s;q) a(s;q)=I$ for some $\gamma$, and
\ref{el-teo-q-lam} -- stating that the eigenvalues of the difference
equation \refe{eqdif} in this case should be of the form
$\lambda_n(q) = C_1\,q^n + C_3$ or $\lambda_n(q)= C_2\,q^{-n}+C_3$.
In section 4 several relevant examples of particular $q$-families of
orthogonal polynomials are illustrated.


\section{Factorization operators}\label{fac-op}

Let introduce a set of functions ${\Phi_n}$
\begin{equation}
\Phi_n(s;q)=d_n^{-1} A(s)\, \sqrt{\rho(s)}\, P_n(s;q),
\label{nor-fun}
\end{equation}
where $d_n$ is the norm of the $q$-polynomials $P_n(s;q)$, $\rho(s)$
is the solution of the Pearson equation \refe{pearson} and $A(s)$
is an arbitrary  continuous function, $A(s)\not=0$ in the interval
$(a,b)$ of orthogonality of $P_n$. If the polynomials $P_n(s;q)$ possess the discrete
orthogonality property  \refe{dis-ort}, then the functions $\Phi_n(s;q)$
satisfy
\bq
\pe{\Phi_n(s;q)}{\Phi_m(s;q)}=\sum_{s=a}^{b-1}{\Phi_n(s;q)}{\Phi_m(s;q)}
\frac{\nabla x_1(s)}{A^2(s)}=\delta_{n,m}.
\label{ort-Phi}
\eq
Notice that if $A(s)=\sqrt{\nabla x_1(s)}$, then the set $(\Phi_n)_n$
is an orthonormal set. Obviously, in the case of a continuous
orthogonality (as for the Askey-Wilson polynomials) one needs to change
the sum in \refe{ort-Phi} by a Riemann integral \cite{nisuuv,atrasu}.\\
Next, we define the $q$-Hamiltonian $\g H(s;q)$ of the form
\begin{equation} \label{gen_hamilt}
\g H(s;q):=
\frac{1}{\nabla x_1(s)}A(s)H(s;q)\frac{1}{A(s)},
\end{equation}
where
\begin{equation}\label{ham3}
\dst H(s;q):=-\frac{\sqrt{\sigma(-s\!-\!\mu\!+\!1)\sigma(s)}}
{\nabla x(s)}e^{-\partial_s}
-\frac{\sqrt{\sigma(-s\!-\!\mu)\sigma(s+1)}}{\Delta x(s)}
e^{\partial_s}+\left(\frac{\sigma(-s\!-\!\mu)}{\Delta
x(s)}+\frac{\sigma(s)}{\nabla x(s)} \right) I,
\end{equation}
$e^{\alpha \partial s} f(s)=f(s+\alpha)$ for all $\alpha \in \CC$
and $I$ is the identity operator. If we now use the identity
$\nabla=\Delta-\nabla \Delta$ and the equation \refe{eqdif}, we find
that
\begin{equation} \label{gen_ham2}
\g H(s;q)\Phi_n(s;q) =\lambda_n \Phi_n(s;q),
\end{equation}
i.e., the functions $\Phi_n(s;q)$, defined in \refe{nor-fun}, are the
eigenfunctions of the associated operator $\g H(s;q)$.
\\
Our first step is to find two operators $a(s;q)$ and $b(s;q)$
such that the Hamiltonian $\g H(s;q) = b(s;q)\,a(s;q)$, i.e.,
the operators  $a(s;q)$ and $b(s;q)$ {\em factorize} the
$q$-Hamiltonian $\g H(s;q)$. But before exhibiting their explicit
form let us point out that if there exists a pair of such operators,
then there are infinitely many of them. Indeed, let $a(s;q)$ and
$b(s;q)$ be such operators that $\g H(s;q)= b(s;q)\,a(s;q)$
and let $U(s;q)$ be an arbitrary unitary operator, i.e., $U^{\dag}(s;q)\,U(s;q)=I$.
Then the operators
$$
\widetilde a(s;q):= U(s;q)\,a(s;q) \, , \quad \quad
\widetilde b(s;q):= b(s;q)\,U^{\dag}(s;q)\, ,
$$
also factorize $\g H(s;q)$ for
$$
\widetilde b(s;q)\,\widetilde a(s;q)= b(s;q)\,U^{\dag}(s;q)
\,U(s;q)\,a(s;q) = b(s;q)\,I\,a(s;q)= b(s;q)\,a(s;q)= \g H(s;q).
$$
This arbitrariness in picking up a particular unitary operator $U(s)$
is very essential because it enables one to construct a closed algebra,
which contains a Hamiltonian $\g H(s;q)$ itself. An explicit form of
the spectrum of this Hamiltonian may then be found by purely algebraic
arguments from the knowledge of representations of this algebra (which
is therefore referred to as {\em a dynamical algebra}).
\\
If one applies the standard factorization procedure to the equation
\refe{eqdif}, then the following difference operators emerge
\bd
Let $\alpha$ be a real number and $A(s)$ an arbitrary
continuous non-vanishing function in $(a,b)$. We define a family of $\alpha$-down
and $\alpha$-up operators by
\begin{equation} \label{alp_oper}
\begin{array}{l}
\!\! \dst \g a^{\downarrow}_\alpha(s;q)\!:=\!\frac{A(s)}{\sqrt{\nabla x_1(s)}}
\, e^{-\alpha \partial_s}
\!\left( e^{\partial_s}\sqrt{\frac{\sigma(s)}{\nabla x(s)}} -
\sqrt{\frac{\sigma(-s-\mu)}{\Delta x(s)}}
\right)\!\frac{1}{A(s)},
 \\[0.6cm]
\!\! \dst \g a^{\uparrow}_\alpha(s;q)\!:=\!\frac{1}{\nabla x_1(s)} A(s)
\!\left( \sqrt{\frac{\sigma(s)}{\nabla x(s)}} e^{-\partial_s} -
\sqrt{\frac{\sigma(-s-\mu)}{\Delta x(s)}}  \right) \!
e^{\alpha \partial_s} \frac{\sqrt{\nabla x_1(s)}}{A(s)},
\end{array}
\end{equation}
respectively.
\ed
A straightforward calculation (by using the simple identity
$e^{\partial_s}\,\nabla = \Delta$) shows that for all
$\alpha\in\RR$
\bq \label{fac}
\g H(s;q)=\g a^{\uparrow}_\alpha(s;q)\g a^{\downarrow}_\alpha(s;q),
\eq
i.e., the operators $\g a^{\downarrow}_\alpha(s;q)$ and
$\g a^{\uparrow}_\alpha(s;q)$ factorize the Hamiltonian, defined
in (\ref{gen_hamilt}). i.e., we have the following
\bt \label{teo-fac}
Given a $q$-Hamiltonian \refe{gen_hamilt} $\g H(s;q)$, then
the operators $\g a^{\downarrow}_\alpha(s;q)$ and
$\g a^{\uparrow}_\alpha(s;q)$ defined in  \refe{alp_oper}
are such that for all $\alpha\in\CC$,
$\g H(s;q)=\g a^{\uparrow}_\alpha(s;q)\g a^{\downarrow}_\alpha(s;q)$.
\et
Our next step is to find a dynamical algebra, associated with the
Hamiltonian $\g H(s;q)$. To this end we will need the following
definition
\bd
A function $f(z)$ is said to be a linear-type function of $z$, if
there exist two functions $F$ and $G$, such that for all $z,\zeta\in\CC$,
this function $f(z)$ can be represented as
$$
f(z+\zeta)=F(\zeta)f(z)+G(\zeta).
$$
\ed
A particular case of the linear-type functions are the $q$-linear
functions, i.e., the functions of the form $f(z)=Aq^z+B$.
For these functions $f(z+\zeta)=F(\zeta)f(z)+G(\zeta)$, where
$F(\zeta)=q^\zeta$ and $G(\zeta)=B(1-q^\zeta)$.
\bn 
If we use the expression in \refe{lam1} for the eigenvalues $\lambda_n$,
then it is straightforward to see that $\lambda_n$ is a $q$-linear function
of $n$ if and only if $\widetilde{\sigma}''=\pm k_q \widetilde{\tau}'$.
Moreover, in this case we have
\begin{equation}\label{lam-lin-rem}
\widetilde{\sigma}''=k_q \widetilde{\tau}' \Rightarrow
\dst\lambda_n(q,+)=\frac{\widetilde{\tau}'}{1-q}(q^n-1),\quad\textrm{and}
\quad \widetilde{\sigma}''=-k_q \widetilde{\tau}' \Rightarrow
\dst \lambda_n(q,-)=\frac{\widetilde{\tau}'}{1-q^{-1}}(q^{-n}-1).
\end{equation}
Notice that $\lambda_n(q,-)=\lambda_n(q^{-1},+)$, i.e., the second
case can be obtained form the first one just by changing $q$ to
$q^{-1}$.
\en
\bp\label{teo1} The function $\lambda_n$ is a $q$-linear function
of $n$ if and only if  it satisfies $\lambda_{n+1}=q\lambda_n+C$.
\ep
\begin{proof}
A straightforward computations show that if $\lambda_n$ is a
$q$-linear function of $n$, then it satisfies the recurrence formula
$\lambda_{n+1}= q\lambda_n+C$, where $C$ is a constant (in this case
$C=\lambda_1$). But the general solution of the difference equation
$\lambda_{n+1}=q\lambda_n+C$ is $\lambda_{n}=Aq^n+D$, where $A$ and $D$
are, in general, non-vanishing constants.
\end{proof}
\bn
Notice that if $\lambda_n$ is a $q$-linear function of $n$, then $\lambda_n$
satisfies the recurrence relation $\lambda_{n+\gamma}-q^{\gamma}\lambda_n=C$
for any  numbers $\gamma$ and $C$.
\en
Finally, we have the following straightforward lemma
\bl
Let $x(s)$ be a $q$-linear function of $s$ and $\lambda_n$ be the
eigenvalue of the difference equation of hypergeometric type \refe{eqdif}.
Then $\lambda_n$ is a $q$-linear function of $n$ if and only if
$\Delta^{(2)}(\sigma(s))=0$ and $q^{-1}$-linear function of $n$ if and
only if $\Delta^{(2)}(\sigma(-s-\mu))=0$, where $\Delta^{(2)}$ is the operator
$\Delta^{(2)}=\frac{\Delta}{\Delta x_1(s)}
\frac{\Delta}{\Delta x(s)}$.
\el
\noindent{\bf Proof:} It follows from equation \refe{lam-lin-rem} and the fact that
$\Delta^{(2)}(\sigma(s))=\frac{[2]_q}{2}(\widetilde{\sigma}''-
\tilde{\tau}' k_q)$ and $\Delta^{(2)}(\sigma(-s-\mu))=
\frac{[2]_q}{2}(\widetilde{\sigma}''+\tilde{\tau}' k_q).$
\hfill $\Box$
\section{Dynamical algebra}\label{alg-din}
We begin this section with the following definition.
\bd
Let $\varsigma$ be a real number, and let $a(s;q)$ and $b(s;q)$ be
two operators. We define the $\mathbf{\varsigma}${\em -commutator
of $a$ and $b$} as
$$
[a(s;q),b(s;q)]_{\varsigma} = a(s;q)b(s;q)-\varsigma b(s;q) a(s;q).
$$
\ed
\bp\label{pro1}
Let $\g H(s;q)$ be an operator, such that there exist two operators
$a(s;q)$ and $b(s;q)$ and two real numbers $\varsigma$ and $\Lambda$,
such that $\g H(s;q)=b (s) a(s;q)$, and $[a(s;q),b(s;q)]_\varsigma=\Lambda$.
Then, if $\Phi(s;q)$ is an eigenvector of the Hamiltonian $\g H(s;q)$,
associated with the eigenvalue $\lambda$, we have
\begin{enumerate}
\item  $\g H(s;q)\{a(s;q)\Phi(s;q)\}=\varsigma^{-1}(\lambda-\Lambda)
\,\{a(s;q)\Phi(s;q)\}$, i.e.,
$a(s;q) \Phi(s;q)$ is the eigenvector of $\g H(s;q)$, associated with the
eigenvalue $\varsigma^{-1}(\lambda-\Lambda)$,
\item  $\g H(s;q)\{b(s;q)\Phi(s;q)\}=(\Lambda+\varsigma\lambda)
\{b(s;q)\Phi(s;q)\}$, i.e.,
$b(s;q)\Phi(s;q)$ is the eigenvector of $\g H(s;q)$, associated with the
eigenvalue $\Lambda+\varsigma\lambda$.
\end{enumerate}
\ep
\noindent{\bf Proof:}
In the first case, since $\g H(s;q)\Phi(s;q) = \lambda \Phi(s;q)$,
$$\ba{rl}
\g H(s;q) \{a(s;q) \Phi(s;q)\}= \!\!\!&b(s;q)a(s;q)\{a(s;q)\Phi(s;q)\}=
\varsigma^{-1}(a(s;q)b(s;q)-\Lambda)
\{a(s;q)\Phi(s;q)\}\\[3mm]
= \!\!\!&\varsigma^{-1}(\lambda-\Lambda)\{a(s;q) \Phi(s;q)\}.
\ea
$$
By the same token, in the second  case
$$\ba{rl}
\g H(s;q) \{b(s;q) \Phi(s;q)\} = \!\!\!& b(s;q)a(s;q)b(s;q) \Phi(s;q)=
b(s;q)(\Lambda + \varsigma\lambda)\Phi(s;q)\\[3mm] = \!\!\!&
(\Lambda +\varsigma\lambda)\{b(s;q) \Phi(s;q)\}.
\ea
$$

\vspace*{-.75cm}

\hfill $\Box$
\\
In the same way one can prove that \bq\label{abphi}
a(s;q)b(s;q)\Phi(s;q)=(\Lambda+\varsigma\lambda)\Phi(s;q). \eq
Moreover, if $\Phi(s;q)$ is an eigenvector of the Hamiltonian $\g
H(s;q)$ (or of the operator $a(s;q)b(s;q)$), then $a^k(s;q)
\Phi(s;q)$ and $b^k(s;q)\Phi(s;q)$ are, in general, also
eigenvectors.
\bn
Obviously, the condition $[a(s;q),b(s;q)]_\q=I$ can be changed to
$[a(s;q),b(s;q)]_\q=\Lambda$, where $\Lambda$ is an arbitrary non-zero
constant. In fact, if the operators $a(s;q)$ and $b(s;q)$
satisfy the $q$-commutation relation $[a(s;q),b(s;q)]_\q=\Lambda$,
then the operators $\g a(s;q)=\Lambda^{-1/2}a(s;q)$ and
$\g b(s;q)=\Lambda^{-1/2}b(s;q)$ satisfy $[\g a(s;q),\g b(s;q)]_\q=I$,
and $\g H(s;q)=\Lambda \g b(s;q)\g a(s;q)$.
\en
The Proposition \ref{pro1} thus refers to the case of a system, described by
a Hamiltonian $\g H(s;q)$, which admits the factorization \refe{fac}
in terms of the operators $a(s;q)$ and $b(s;q)$, satisfying
the $q$-commutation relation $[a(s;q),b(s;q)]_\q = I$. Moreover,
it tells us how to construct a dynamical symmetry algebra for such a
case in a direct fashion \cite{kuda}.  Indeed, let us assume that
$[a(s;q),b(s;q)]_\varsigma = I, \, \varsigma=q^2 \,(or \,q^{-2})$,
and $b(s;q)= a^\dag(s;q)$. Then one can rewrite the $q^2$-commutator
$a(s;q)\,a^\dag(s;q)-q^2\,a^\dag(s;q)\,a(s;q)=I$ in the following form
$$
[a(s;q),a^\dag(s;q)] := a(s;q)\,a^\dag(s;q)-a^\dag(s;q)\,a(s;q)
=I - (1-q^2)a^\dag(s;q)\,a(s;q) := q^{2N(s)}, 
$$
where, by definition, the operator $N(s)$ is equal to
$N(s)=\ln[I-(1-q^2)\,a^\dag(s;q)\,a(s;q)]/\ln q^2$. From this definition
of $N(s)$ it follows that
\bq
[N(s),a(s;q)]=-a(s;q), \quad [N(s),a^\dag(s;q)]=a^\dag(s;q), \label{(b)}
\eq
i.e., $N(s)$ is the number operator. The next (and final) step is to
introduce a new set of the operators
$$b(s;q):=q^{-N(s)/2}\,a(s;q),\quad b^\dag(s;q):=a^\dag(s;q)\,q^{-N(s)/2},
$$
which satisfy the following commutation relation
$$ b(s;q)\,b^\dag(s;q)-q\,b^\dag(s;q)\,b(s;q)=q^{-N(s)}, 
$$
readily derived with the aid of \refe{(b)}. The operators
$b(s;q)$, $b^\dag(s;q)$, and $N(s)$ directly lead to the dynamical
algebra $su_q(1,1)$ with the generators
$$ K_0(s)=\frac{1}{2}\,[N(s)+1/2], \quad K_+(s)=\beta \, (b^\dag(s;q))^2,
\quad K_-(s)=\beta \, b^2(s;q), \quad \beta^{-1}=q+q^{-1}.
$$
It is straightforward to verify that thus defined generators satisfy
the standard commutation relations
$$ [K_0(s),K_{\pm}(s)]= \,\pm \, K_{\pm}(s), \quad
[K_-(s),K_+(s)]= [2K_0(s)]_{q^2},
$$
of the algebra $su_q(1,1)$ (see e.g. \cite{vikl}).

Thus in the case when the operators $a(s;q)$ and $b(s;q)$,
which factorize the Hamiltonian $\g H(s;q)$, satisfy the $q$-commutation
relation $[a(s;q),b(s;q)]_{\q^2} = I$ and $b(s;q)=a^\dag(s;q)$, the appropriate dynamical
symmetry algebra is $su_\q(1,1)$. So the question arises: what are
conditions for inisuuvring that such $q$-commutator takes place? In other words,
we have the following\\
\noindent\textbf{Problem 1:} {\em To find two operators $a(s;q)$ and
$b(s;q)$ and a constant $\q$ such that the Hamiltonian $\g
H(s;q)=b(s;q)a(s;q)$ and
$[a(s;q),b(s;q)]_{\q}=I$.}\\
For the first part we already have the answer (see Theorem
\ref{teo-fac}). The solution of the second one is formulated in the
following two theorems.
\bt \label{el-teo-q-lin} Let  $\g H(s;q)$ be
the following difference operator ($q$-Hamiltonian)
\begin{equation} \label{ham_lin}
\g H(s;q)= \frac{1}{\nabla
x_1(s)}A(s)H(s;q)\frac{1}{A(s)}.
\end{equation}
The operators $b(s;q)=\g a^{\uparrow}_\alpha(s;q)$ and $a(s;q)=\g
a^{\downarrow}_\alpha(s;q)$ given in  \refe{alp_oper}
factorize the Hamiltonian $\g H(s;q)$ \refe{ham_lin} and satisfy
the commutation relation $[a(s;q),b(s;q)]_{\q}=\Lambda$ for a certain
real number $\q$ if and only if the following two conditions hold:
\begin{equation} \label{first_cond_lin}
\frac{\nabla x(s)}{\nabla x_1(s-\alpha)} \sqrt{\frac{\nabla
x_1(s-1)\nabla x_1(s)}{\nabla x(s-\alpha)\Delta x(s-\alpha)}}
\sqrt{\frac{\sigma(s-\alpha)\sigma(-s-\mu+\alpha)}{\sigma(s)\sigma(-s-\mu+1)}}=\q,
\end{equation}
and
\begin{equation} \label{second_cond_lin}
\frac{1}{\Delta
x(s-\alpha)}\!\left(\frac{\sigma(s-\alpha+1)}{\nabla
x_1(s-\alpha+1)}+ \frac{\sigma(-s-\mu+\alpha)}{\nabla
x_1(s-\alpha)}\right)\!-\varsigma\,\frac{1}{\nabla
x_1(s)}\left(\!\frac{\sigma(s)}{\nabla
x(s)}+\frac{\sigma(-s-\mu)}{\Delta x(s)}\!\right)=\Lambda.
\end{equation}
\et \noindent{\bf Proof:} Taking the expression of the operators
$a_\alpha^\uparrow(s)$ and $a_\alpha^\downarrow(s)$, a
straightforward calculus shows that
$a_\alpha^\downarrow(s)a_\alpha^\uparrow(s)=
A_1(s)e^{\partial_s}+A_2(s)e^{-\partial_s} +A_3(s)I$, where
\bq\label{help1}
\begin{array}{l}
\dst A_1(s)=-\sqrt{\frac{\nabla x_1(s+1)}{\nabla x_1(s)}}\frac{A(s)}{A(s+1)}
\sqrt{\frac{\sigma(s+1-\alpha)\sigma(-s-\mu-1+\alpha)}
{\Delta x(s-\alpha)\Delta x(s+1-\alpha)}}\frac{1}{\nabla x_1(s+1-\alpha)},\\[0.5cm]
\dst A_2(s)=-\sqrt{\frac{\nabla x_1(s-1)}{\nabla x_1(s)}}\frac{A(s)}{A(s-1)}
\sqrt{\frac{\sigma(s-\alpha)\sigma(-s-\mu+\alpha)}
{\Delta x(s-1-\alpha)\Delta x(s-\alpha)}}\frac{1}{\nabla x_1(s-\alpha)}, \\[0.5cm]
\dst A_3(s)=\frac{1}{\Delta
x(s-\alpha)}\left[\frac{\sigma(s+1-\alpha)} {\nabla
x_1(s+1-\alpha)}+\frac{\sigma(-s-\mu+\alpha)}{\nabla
x_1(s-\alpha)}\right].
\end{array}
\eq
In the same way, using \refe{ham3} and \refe{gen_hamilt}
we have $a_\alpha^\uparrow(s)a_\alpha^\downarrow(s)=\g H(s;q)=
B_1(s)e^{\partial_s}+B_2(s)e^{-\partial_s}+B_3(s)I$, where
\bq\label{help2}
\begin{array}{l}
\dst B_1(s)=-\frac{1}{\nabla x_1(s)}\frac{A(s)}{A(s+1)}
\frac{\sqrt{\sigma(-s-\mu)\sigma(s+1)}} {\nabla x(s+1)},\\[0.5cm]
\dst B_2(s)=-\frac{1}{\nabla x_1(s)}\frac{A(s)}{A(s-1)}
\frac{\sqrt{\sigma(-s-\mu+1)\sigma(s)}}
{\nabla x(s)},\\[0.5cm]
\dst B_3(s)=\frac{1}{\nabla x_1(s)}\left[\frac{\sigma(s)}{\nabla
x(s)}+ \frac{\sigma(-s-\mu)}{\Delta x(s)}\right].
\end{array}
\eq
Consequently,
\bq\label{help3}
[a_\alpha^\downarrow(s),a_\alpha^\uparrow(s)]_\varsigma=
\Big(A_1(s)-\varsigma B_1(s)\Big)e^{\partial_s}+\Big(A_2(s)-
\varsigma B_2(s)\Big)e^{-\partial_s} +\Big(A_3(s)-\varsigma
B_3(s)\Big)I.
\eq
To eliminate the two terms in the right-hand side of \refe{help3},
which are proportional to the difference operators $\exp(\pm \partial_s)$,
one must require that
\bq
A_1(s) - \varsigma\,B_1(s)=0, \quad \quad
A_2(s)- \varsigma\,B_2(s)=0.             \label{(25)}
\eq
Off hand, it is not evident that one can satisfy both of the
relations \refe{(25)}, which only involve the same constant $\varsigma$.
But it is straightforward to verify from \refe{help1} and \refe{help2}
that
$$ A_1(s)\,B_2(s+1)=A_2(s+1)\,B_1(s),                     
$$
or, equivalently,
$$ \frac{A_1(s)}{B_1(s)}=\frac{A_2(s+1)}{B_2(s+1)}.
$$
Hence, the requirement that $A_1(s)=\varsigma\,B_1(s)$ entails the
relation $A_2(s)=\varsigma\,B_2(s)$, and vice versa.
{From} \refe{help3} it is now evident that the commutator
$[a_\alpha^\uparrow(s),a_\alpha^\downarrow(s)]_\q$ is a constant
if \refe{(25)} holds and the factor $ A_3(s)-\varsigma \,B_3(s)$ is a
constant. Thus, the required conditions \refe{first_cond_lin} and
\refe{second_cond_lin} immediately follow. \hfill $\Box$

\bt\label{el-teo-q-lam} Let $(\Phi_n)_n$ be the eigenfunctions
of  $\g H(s;q)$ corresponding to the eigenvalues $(\lambda_n)_n$
and suppose that the problem 1 has a solution for $\Lambda\neq0$.
Then, the eigenvalues $\lambda_n$ of the difference equation
\refe{gen_ham2} are $q$-linear or
$q^{-1}$-linear functions of $n$, i.e., $\lambda_n=
C_1 q^n+C_3$ or $\lambda_n=C_2q^{-n}+C_3$, respectively.
\et
\bdm\footnote{For an alternative proof in the case $\alpha=0$ see the
appendix.} In the following we use the notation
$\q=q^\gamma$. Suppose that problem 1 has a solution with
 $\Lambda\neq0$ and $\lambda_n$ is not a $q$-linear
(respectively, $q^{-1}$) function of $n$. From Proposition
\ref{pro1} we know that
$\g a^{\uparrow}_\alpha(s;q)\Phi_n(s;q)$ is and eigenvector
of $\g H(s;q)$ corresponding to the eigenvalue
$\Lambda+q^\gamma \lambda_n$. If we denote by $\Phi_{m(n);q}$ such
eigenvector where $m(n)$ is a function of $n$, then we have
$\Lambda+q^\gamma\lambda_n=\lambda_{m(n)}$. Then using
\refe{lam1} we get
\bq\label{lam-h1}
\lambda_{m(n)}= C_1 q^{m(n)}+C_2q^{-m(n)}+C_3,\quad
C_1C_2=L_q.
\eq
On the other hand
\bq\label{lam-h2}
\lambda_{m(n)}=\Lambda+q^\gamma \lambda_n=
C_1 q^{\gamma}q^n+C_2q^{\gamma}q^{-n}+q^\gamma C_3+\Lambda=
C_1'q^n+C_2'q^{-n}+C_3'.
\eq
But here, since $\lambda_{m(n)}$ is an eigenvalue of \refe{eqdif},
again we have the condition $C_1'C_2'=L_q$, thus
$C_1C_2=L_q=C_1'C_2'=C_1C_2q^{2\gamma}$  so $q^{2\gamma}=1$, i.e.,
$\gamma=0$, or $C_1C_2=0$. In the first case, equating \refe{lam-h1} and
\refe{lam-h2}, we have that
$C_3'=C_3q^\gamma+\Lambda=C_3$, i.e., $\Lambda=0$ that is a
contradiction. Thus $C_1C_2=0$ from where the result easily
follows.
\edm
\\
It is worth noting that the $q$-linearity of the eigenvalues
of $\g H(s;q)$ is only the necessary condition, i.e., it is not
sufficient. So there are cases when the eigenvalues $\lambda_n(q)$
are $q$-linear (for instance, those which correspond to the
$q$-Meixner, the $q$-Charlier, and the $q$-Laguerre polynomials
with $a\neq q^{-1/2}$), but the corresponding $q$-Hamiltonians
$\g H(s;q)$ do not admit the factorization in terms of $q$-commuting
operators. This just
reflects the fact that an appropriate dynamical algebra is not
$su_q(1,1)$ and one has to consider a more complicated quadratic
algebra $AW(3)$ \cite{Zed}. The problem of finding an explicit
connection between the generators of the algebra $AW(3)$ and the
operators $a(s;q)$ and $b(s;q)$, which factorize the
$q$-Hamiltonians for these cases, will be attended in a separate
publication.
\bn\label{rem-red-no-lin} A special important case of the
non-linear lattice is when $x(s)=\half(q^s+q^{-s})$. In this case
if we put $\alpha=\half$, then the conditions \refe{first_cond_lin} and
\refe{second_cond_lin} of Theorem \ref{el-teo-q-lin} becomes
$$
\sqrt{\frac{\sigma(s-\half)\sigma(-s+\half)}
{\sigma(s)\sigma(-s+1)}}=\q,
$$
and
$$
\frac{1}{\nabla x_1(s)}\!\left(\frac{\sigma(s+\half)}{\Delta
x(s)}+ \frac{\sigma(-s+\half)}{\nabla x(s)}\right)\!-\!\varsigma
\frac{1}{\nabla x_1(s)}\left(\!\frac{\sigma(s)}{\nabla
x(s)}+\frac{\sigma(-s)}{\Delta x(s)}\!\right)=\Lambda.
$$
respectively. Moreover, if we put $A(s)=1$ then,
$\g H(s;q)=(\nabla x_1(s))^{-1}H(s;q)$ and the $\alpha$-operators
simplify
$$\ba{l}\dst
a_{1/2}^\downarrow(s)=\frac{1}{\nabla x_1(s)}
\left(e^{\frac12\partial_s}\sqrt{\sigma(s)}-
e^{-\frac12\partial_s}\sqrt{\sigma(-s)}\right),
\\[5mm]\dst
a_{1/2}^\uparrow(s)=\frac{1}{\nabla x_1(s)}
\left(\sqrt{\sigma(s)} e^{-\frac12\partial_s}-
\sqrt{\sigma(-s)}e^{\frac12\partial_s}\right).
\ea
$$
\en
Now we can formulate the
\\[0.6cm]
\noindent\textbf{Problem 2:}
{\em To find two operators $a(s;q)$ and $b(s;q)$ and a constant $\q$
such that the Hamiltonian $\g H(s;q)=b(s;q)a(s;q)$ and
$[a(s;q),b(s;q)]_{\q}=I$ and such that $a(s;q)$ and
$b(s;q)$ are the lowering and raising operators, i.e.,}
\bq
a(s;q) \Phi_n(s;q)=D_n \Phi_{n-1}(s;q) \quad \mbox{and}
\quad b(s;q)\Phi_n(s;q)=U_n \Phi_{n+1}(s;q).
\label{low-rai1}
\eq
Again, without loss of generality, we will change the condition
$[a(s;q),b(s;q)]_{\q}=I$ into
$$
[a(s;q),b(s;q)]_{\q}=\Lambda \quad \mbox{and chose} \quad \Lambda=\lambda_1.
$$
Also the operators
$b(s;q)=\g a^{\uparrow}_\alpha(s;q)$ and
$a(s;q)=\g a^{\downarrow}_\alpha(s;q)$, given in \refe{alp_oper},
provide the factorization of $\g H(s;q)$.
\\
If we now apply $b(s;q)$ to the first equation of \refe{low-rai1}
and use the second one as well as \refe{gen_ham2}, we find that
$\lambda_n=D_n U_{n-1}$. On the other hand, applying $a(s;q)$
to the second equation in \refe{low-rai1} and using the first
one as well as \refe{abphi} we obtain
$\lambda_1+\q\lambda_n=U_nD_{n+1}=\lambda_{n+1}$. Thus, using
Proposition \ref{teo1} we conclude that $\lambda_n$ should be a
$\q$-linear function, i.e., $\lambda_n$ has the form
$\lambda_n=A\q^n+D$, where $A$ and $D$ are non-vanishing constants.
Moreover, using the recurrence $\lambda_1+\q\lambda_n=\lambda_{n+1}$
and Proposition \ref{pro1} we obtain that
$$
\ba{l}
\g H(s;q)\{a(s;q)\Phi_n(s;q)\}=\lambda_{n-1}\{a(s;q)\Phi_n(s;q)\},\\[0.3cm]
\g H(s;q)\{b(s;q)\Phi_n(s;q)\}=\lambda_{n+1}\{b(s;q)\Phi_n(s;q)\},
\ea
$$
i.e., $a(s;q)\Phi_n(s;q)$ is the eigenvector associated to the
eigenvalue $\lambda_{n-1}$ and $b(s;q)\Phi_n(s;q)$ is the eigenvector
associated to the eigenvalue $\lambda_{n+1}$, so they really
are the lowering and raising operators associated to the
Hamiltonian $\g H(s;q)$.
\\
The next important question to be considered is when the
$\alpha$-operators are mutually adjoint? Obviously the answer to
this question depends on how the scalar product is defined. As an
example, let us consider the discrete orthogonality case
\refe{ort-Phi} and compute  $\pe{\g a^{\downarrow}_\alpha
\Phi_{n+1}}{\Phi_k}$:
\begin{small}
\[
\begin{split} &
\pe{\g a^{\downarrow}_\alpha(s;q) \Phi_{n+1}(s;q)}{\Phi_k(s;q)}=\\[5mm]
=& {\dst\sum_{s=a}^{b-1}}
\frac{\sqrt{\nabla x_1(s)\rho(s)}}{d_n d_{n+1}}
\!\left[e^{-\alpha \partial_s}\!\left(
\sqrt{\frac{\sigma(s+1)\rho(s+1)}{\nabla x(s+1)}}
P_{n+1}(s\!+\!1;q)\!-\!
\sqrt{\frac{\sigma(-s-\mu)\rho(s)}{\Delta x(s)}}P_{n+1}(s;q)
\!\right)\!\right]\!P_{k}(s;q)\\[5mm]
=& \dst\sum_{s=a}^{b-1}\frac{\sqrt{\nabla x_1(s)\rho(s)}}{d_n d_{n+1}}
\left[e^{-\alpha \partial_s}
\!\left( \sqrt{\frac{\sigma(-s-\mu)\rho(s)}{\Delta x(s)}}
\Delta P_{n+1}(s;q)\right)\right]
P_{k}(s;q)\\[5mm]
=& \dst\sum_{s=a}^{b-1}
\frac{\sqrt{\nabla x_1(s)\rho(s)}}{d_n d_{n+1}}
\sqrt{\frac{\sigma(-s-\mu+\alpha)\rho(s-\alpha)}
{\Delta x(s-\alpha)}}\Delta P_{n+1}(s-\alpha;q) P_{k}(s;q),
\end{split}
\]
\end{small}
where the second equality follows from the Pearson equation
\refe{pearson}. Next we have
\begin{small}
\[
\begin{split}
\pe{\Phi_{n+1}(s;q)}{\g a^{\uparrow}_\alpha (s)\Phi_{k}(s;q)}=&
\dst\underbrace{\sum_{s=a}^{b-1}
\frac{\sqrt{\sigma(s)\rho(s)}}{d_n d_{n+1}\sqrt{\nabla x(s)}}
\!\left[e^{-\partial_s}
\left(e^{\alpha \partial_s}\sqrt{\nabla x_1(s)\rho(s)}
P_{k}(s;q)\right)\right]P_{n+1}(s;q)}_{S_1}\\[5mm] &
\dst-\underbrace{\sum_{s=a}^{b-1}
\frac{\sqrt{\sigma(-s-\mu)\rho(s)}}{d_n d_{n+1}\sqrt{\Delta x(s)}}
\left[e^{\alpha\partial_s}\left(\sqrt{\nabla x_1(s)\rho(s)}
P_{k}(s;q)\!\right)\!\right]\!P_{n+1}(s;q)}_{S_2}
\end{split}
\]
\end{small}
If we use the boundary conditions $\sigma(s) \rho(s)\Big|_{s=a,b}
= 0$, then the sum $S_1$ becomes
\[
\begin{split}
S_1=& \dst \sum_{s=a+1}^{b}
\frac{\sqrt{\sigma(s)\rho(s)}}{d_n d_{n+1}\sqrt{\nabla x(s)}}
\!\left[e^{-\partial_s}
\left(e^{\alpha \partial_s}\sqrt{\nabla x_1(s)\rho(s)}
P_{k}(s;q)\right)\right]P_{n+1}(s;q) \qquad (s\to s+1)\\[5mm]
=& \dst \sum_{s=a}^{b-1}
\frac{\sqrt{\sigma(s+1)\rho(s+1)}}{d_n d_{n+1}\sqrt{\nabla x(s+1)}}
\!\left[\left(e^{\alpha \partial_s}\sqrt{\nabla x_1(s)\rho(s)}
P_{k}(s;q)\right)\right]P_{n+1}(s+1;q) \\[5mm]
=& \dst \sum_{s=a}^{b-1}
\frac{\sqrt{\sigma(-s-\mu)\rho(s)}}{d_n d_{n+1}\sqrt{\Delta x(s)}}
\!\left[e^{\alpha \partial_s}\left(\sqrt{\nabla x_1(s)\rho(s)}
P_{k}(s;q)\right)\right]P_{n+1}(s+1;q),
\end{split}
\]
where the last equality holds due to the Pearson equation
\refe{pearson}. Now subtracting $S_1-S_2$ yields
\begin{small}
$$\ba{rl}
\pe{\Phi_{n+1}(s;q)}{\g a^{\uparrow}_\alpha (s)\Phi_{k}(s;q)}& \!\!\!=\dst
 \dst \sum_{s=a}^{b-1} \frac{\sqrt{\sigma(-s-\mu)\rho(s)}}{d_n d_{n+1}\sqrt{\Delta x(s)}}
\!\left[e^{\alpha \partial_s}\left(\sqrt{\nabla x_1(s)\rho(s)}
 P_{k}(s;q)\right)\right] \Delta P_{n+1}(s;q) \\[5mm]
=&\!\!\! \dst\sum_{s=a}^{b-1}
\frac{\sqrt{\sigma(-s-\mu)\rho(s)}}{d_n d_{n+1}\sqrt{\Delta x(s)}}
\sqrt{\nabla x_1(s+\alpha)\rho(s+\alpha)}
 P_{k}(s+\alpha;q) \, \Delta P_{n+1}(s;q).
\ea
$$
\end{small}
Then, in the discrete case, a sufficient condition for the operators
$\g a^{\uparrow}_\alpha(s;q)$ and $\g a^{\downarrow}_\alpha(s;q)$
to be mutually adjoint, i.e., $\pe{\g a^{\downarrow}_\alpha
(s)\Phi_{n+1}(s;q)}{\Phi_k(s;q)}=\pe{\Phi_{n+1}(s;q)}
{\g a^{\uparrow}_\alpha (s)\Phi_{k}(s;q)}$, is that $\alpha=0$.
For the general case as well as for the discrete cases when
$\alpha\neq0$, the problem requires a more detailed study of each case.
A discussion of these cases will be considered elsewhere.

\section{Examples}

In this section we will exhibit several examples involving
some well-known families of $q$-po\-ly\-no\-mials. Let us point out
that from the whole $q$-Askey tableau \cite{kost} we need only to
consider those cases, when $\lambda_n$ is a $q$-linear or $q^{-1}$-linear
function of $n$ (see Theorem \ref{el-teo-q-lam}).

\subsection{The linear lattice $x(s)=c_1 q^s+c_3$: The $q$-Hahn
Tableau}
We start with the so-called $q$-Hahn tableau
(see e.g. \cite{alme,kost} and references therein).
Taking into account \ref{el-teo-q-lin}, the $\alpha$-operators are chosen in
such a way that $A(s)=\sqrt{\nabla x_1(s)}$, so the orthonormal
functions are defined by
$$
\Phi_n(s;q)={d_n}^{-1}\sqrt{\rho(s)\nabla x_1(s)}P_n(s;q)
$$
and satisfy the discrete orthogonality relation
$$
 \sum_{s_i=a}^{b-1} \Phi_n(s_i;q) \Phi_m(s_i;q)
=\delta_{nm}.
$$
By using the theorem \ref{el-teo-q-lin}, one can solve the factorization
problem 1 for the $q$-Hahn tableau.

Before starting with some relevant examples, we recall that for
the $q$-linear lattice $x(s)=c_1 q^s+c_3$,
$$
\g e_q(s):=\frac{\nabla x(s)}{\nabla x_1(s-\alpha)}
\sqrt{\frac{\nabla x_1(s-1)\nabla x_1(s)}{\nabla x(s-\alpha)\Delta
x(s-\alpha)}} = q^{2\alpha-1},
$$
whereas if $x(s)=c_1 q^{-s}+c_3$, then $\g e_q(s)=q^{-2\alpha+1}$.

\bn
Let $x(s)=q^s$. Then in the case when  $\sigma(s)$
or  $\sigma(-s-\mu)=\sigma(s)+\tau(s)\nabla x_1(s)$
are constants, the operators \refe{alp_oper}
define a dynamical algebra if and only if
$\alpha =1$ and $\q=q$ or $\alpha=0$ and $\q=q^{-1}$,
respectively. To prove this
assertion it is sufficient to use the formula
(\ref{first_cond_lin}) that yields, for the first case
$$
\g e_q(s)\sqrt{\frac{\sigma(-s-\mu+\alpha)}
{\sigma(-s-\mu+1)}}=
q^{2\alpha-1}\sqrt{\frac{\sigma(-s-\mu+\alpha)}
{\sigma(-s-\mu+1)}}=\q.
$$
Choosing $\alpha=1$ we have that $\q=q$. The other case is
analogous. These cases constitute the most simple ones.
\en

\subsubsection{Stieltjes-Wigert $q$-polynomials $S_n(x;q)$}

The  Stieltjes-Wigert functions in $x(s)=q^s$, i.e. the functions,
associated with the Stieltjes-Wigert polynomials, are defined by
$$
\Phi_n(x;q)= \frac1{d_n}\sqrt{(-x,-q/x;q)_\infty}\
{_1\phi_1}\left( \begin{array}{c} q^{-n} \\
0 \end{array} \Bigg\vert q;-xq^{n+1}\right),\quad
x(s)=q^s,
$$
$$
d_n=\frac{q^{n/2}}{(q;q)_\infty}\sqrt{\frac{(q^{n+1};q)_\infty}{\log q^{-1}}}\ .
$$
In this case we have chosen $A(s)=\sqrt{\nabla x_1(s)}$.
The above functions $\Phi_n(x;q)$ possess the following orthogonality
property $\int_0^\infty\Phi_n(x;q)\Phi_m(x;q)dx=\delta_{n,m}$.

Since for the Stieltjes-Wigert $q$-polynomials $\sigma(s)=q^{s-1}$
and $\sigma(s)+\tau(s)\nabla x_1(s)=q^{2s}$ \cite{kost},
one can define the following Hamiltonian
$$
\g H(s;q)=\frac{1}{(1-q)}\left\{
\left(1+q^{-s}\right)I
-q^{-(s+1)/2} e^{\partial_s} - q^{-s/2} \ e^{-\partial_s}\right\},
$$
for which one has $\g H(s;q)\Phi_n(q^s;q)=\frac{1-q^n}{1-q}\,\Phi_n(q^s;q)$.

Next we check the conditions of the Theorem \ref{el-teo-q-lin}:
The first condition \refe{first_cond_lin} yields $\q=q^{\alpha/2}$
and substituting it in the second one, we obtain that it holds when
$\alpha =2$, i.e.,  $\varsigma = q$ and therefore
$$
\begin{array}{l}
\dst \g a^{\downarrow}(s;q) := \g
a^{\downarrow}_2(s;q)=\frac{1}{\sqrt{1-q}}\left(e^{-2\partial_s}
-e^{-\partial_s} q^{-s/2}\right), \\[0.6cm]
\dst \g a^{\uparrow}(s;q) := \g
a^{\uparrow}_2(s;q)=
\frac{1}{\sqrt{1-q}}\left(e^{2\partial_s}-q^{-s/2} e^{\partial_s} \right),
\end{array}
$$
$\g H(s;q)= \g a^{\uparrow}(s;q) \g a^{\downarrow}(s;q)$ and
$[\g a^{\downarrow}(s;q), \g a^{\uparrow}(s;q)]_{q}=I$.
It is not hard to verify that  in this case the operators $\g a^{\uparrow}(s;q)$ and
$\g a^{\downarrow}(s;q)$   are the lowering and raising operators for the
functions $\Phi_n(q^s;q)$ (see   \cite{kost}, p.117, (3.27.6) and (3.27.8)).
We remind the reader that the moment problem, associated with the Stieltjes-Wigert
polynomials, is indeterminate \cite{akh,shta} and therefore there are several distinct
weight functions (both continuous and discrete ones), with respect to which they are
orthogonal. A similar result for the case of a discrete orthogonality condition has been
first considered in  \cite{atsu2}.

\subsubsection{Al-Salam \& Carlitz I and II  $q$-polynomials $U_n^{(a)}(x;q)$ and $V_n^{(a)}(x;q)$}

The Al-Salam--Carlitz polynomials of type I and of type II,
depend on an additional parameter $a$ and therefore they occupy the
next level in the Askey scheme (see \cite[p.114]{kost}). Since these two
families are interrelated,
$$
V_n^{(a)}(x;q) = U_n^{(a)}(x;q^{-1}),
$$
it is sufficient to consider only one of them.

Let us define the functions \cite[p. 114]{kost}
$$
\Phi_n(s;q)= \frac1{d_n} \frac{a^{s/2}q^{s^2/2}}{\sqrt{(q,aq;q)_s}}
\, _2\phi_0\left(\begin{array}{c|c}q^{-n},q^{-s} \\[-0.35cm]
& q;\dst \frac{q^n}{a} \\[-0.2cm] - \end{array} \right), \quad a>0,
$$
$$
d_n= (-1)^n(aq)^{-n/2}\sqrt{\frac{(q;q)_n}{(aq;q)_\infty}},
$$
which satisfy the discrete orthogonality relation
$\sum_{k=0}^\infty \Phi_n(k;q)\Phi_m(k;q)=\delta_{m,n}$.
Since for this family $\sigma(s)=(q^{-s}-1)(q^{-s}-a)$ and
$\sigma(s)+\tau(s)\nabla x_1(s)=a$, where $x(s)=q^{-s}$, then one
derives that the difference Hamiltonian $\g H(s;q)$ has the following form
$$\ba{l}
\g H(s;q)=\dst \frac{1}{1-q}\,\left[a\,q^{2s+1} + (1-q^s)(1-a\,q^s)\right.\\[4mm]
\qquad \qquad \dst \left.-\, e^{\partial_s}\sqrt{a(1-q^s)(1-a\,q^s)}\,q^{s-1/2}-
\sqrt{a(1-q^s)(1-a\,q^s)}\,q^{s-1/2}\,e^{-\partial_s}\,\right]\,.
\ea
$$
Thus $\g H(s;q)\,\Phi_n(s;q) =\frac{1-q^n}{1-q}\,\Phi_n(s;q)$.
This Hamiltonian is factorized,
$\g H(s;q) =\g a^{\uparrow}(s;q)\g a^{\downarrow}(s;q)$,
in terms of the difference operators
$$
\g a^{\downarrow}(s;q):=\g a^{\downarrow}_0(s;q)=
\frac{1}{\sqrt{1-q}}\,\left[\sqrt{a}\,q^{s+1/2}-
e^{\partial_s}\,\sqrt{(1-q^s)(1-a\,q^s)}\,\right],
$$
$$
\g a^{\uparrow}(s;q):=\g a^{\uparrow}_0(s;q)=
\frac{1}{\sqrt{1-q}}\,\left[\sqrt{a}\,q^{s+1/2}-
\sqrt{(1-q^s)(1-a\,q^s)}\,e^{-\partial_s}\,\right]\,,
$$
that satisfy the commutation relation
$[\g a^{\downarrow}(s;q), \g a^{\uparrow}(s;q)]_{q}=I$.
Thus the dynamical algebra for this family is
also $su_q(1,1)$. As in the previous case, the operators $\g a^{\downarrow}_2(s;q)$ and
$\g a^{\uparrow}(s;q)$   are the lowering and raising operators for the
functions $\Phi_n(s;q)$ (see   \cite{kost}).

A special case of the Al-Salam\&Carlitz $q$-polynomials of type II
are the discrete Hermite $q$-polynomials $\widetilde{h}_n(x;q)=
{\rm{i}}^{-n}\,V^{(-1)}_n({\rm i}\,x;q)$ of type II.

All these cases are closely related to some models of
$q$-oscillators \cite{mac,atsu2,atsu1,assu1,assu2,atfrwo}.

\subsubsection{Wall polynomials $p_n(x;a|q)$}
Our next example is the little $q$-Laguerre / Wall polynomials.
In this case we define the function
$$
\Phi_n(s;a;q)=  \frac1{d_n(a;q)}\frac{(aq)^{s/2}}{\sqrt{(q;q)_s}}\
{_2\phi_1}\left( \begin{array}{c} q^{-n},0 \\
aq \end{array} \Bigg\vert q;q^{s+1}\right),
$$
where $x:=x(s)=q^s$ and $A(s)=\sqrt{\nabla x_1(s)}$.

This function satisfies the discrete orthogonality relation  $\sum_{k=0}^\infty
\Phi_n(k;q)\Phi_m(k;q)=\delta_{n,m}$, provided that
$$
d_n(a;q)=\frac{(aq)^{n/2}}{(aq;q)_\infty} \, \sqrt{(q;q)_n (aq^{n+1};q)_\infty}\ .
$$

In this case  $\sigma(s)=q^{-1}q^s(q^s-1)$ and $\sigma(s)+\tau(s)\nabla x_1(s)=-\,aq^s$,
thus the corresponding Hamiltonian is
$$
\g H(s;a;q)=
\frac{1}{(1-q)x}\left({q}(a+1-x)I-\sqrt{aq(1-qx)}\,e^{\partial_s}-q\sqrt{aq(1 -x)}\,e^{-\partial_s}\right)
$$
and $\g H(s;a;q) \Phi_n(s;a;q)=\frac{1-q^{-n}}{1-q^{-1}}\,\Phi_n(s;a;q)$.

Then, the first condition of Theorem \ref{el-teo-q-lin}
leads to the value $\alpha=0$ and $\q=q^{-1/2}$ and the
second condition holds if and only if the parameter $a$
of the above functions is equal to $q^{-1/2}$,
but they do not lead to the lowering and raising operators.

So to introduce the lowering and raising operators one has to
consider the following operators
$$
\ba{l}
\dst \g a(s;a;q) :=
\frac1{\sqrt{\left(1-q\right)x}} \left(\sqrt{(1-qx)}\, e^{\partial_s} -\sqrt{aq}\,I\right),  \\[0.6cm]
\dst \g a^{\dag}(s;a;q) := \frac1{\sqrt{\left(1-q\right)x}}
\left(\sqrt{q(1-x)}\, e^{-\partial_s}-\sqrt{aq}\,I \right).
\ea
$$
The above mutually adjoint operators factorize the Hamiltonian $\g H(s;a;q)$, i.e.,
 $\g H(s;a;q)=\g a^{\dag}(s;a;q) \, \g a(s;a;q)$, and
 they satisfy the commutation relation
$$
\g a(s;a/q;q)  \g a^{\dag}(s;a/q;q) - q^{-1} \g a^{\dag}(s;a;q)  \g a(s;a;q)=I.
$$
{From} this relation it follows that their action on the functions $\Phi(s;a;q)$
is given by
\bq\label{lr-wall}
\begin{split}
&\g a(s;a;q)\,\Phi_n(s;a;q)=\sqrt{\frac{1-q^{-n}}{1-q^{-1}}}\, \Phi_{n-1}(s;aq;q),
\\[4mm]
&\g a^{\dag}(s;a/q;q)\,\Phi_n(s;a;q)=\sqrt{\frac{1-q^{-n-1}}{1-q^{-1}}}\, \Phi_{n+1}(s;a/q;q).
\end{split}
\eq
To verify these formulae one needs to use the property of the normalization constant $d_n(a;q)$ that
$$
d_n(a;q)=\frac{(1-a)\,q^{n/2}}{\sqrt{a(1-q^{n+1})}}\,d_{n+1}(a/q;q).
$$
Observe that the operators $\g a(s;a;q)$ and $\g a^{\dag}(s;a;q)$ not only lower and raise, respectively,
the index $n$, but they alter also the parameter $a$. That is why an appropriate dynamical algebra
in this case is not $su_q(1,1)$ (see the discussion above, which follows the proof of Theorem
\ref{el-teo-q-lin} in section \ref{alg-din}). The formulae \refe{lr-wall} are equivalent to proving that the
forward and backward shift operators for the little $q$-Laguerre polynomials have the form
(see \cite{kost} p. 107, (3.20.6) and (3.20.8))
$$
(e^{\partial_s}-I)\,p_n(q^s;a;q)=q^{s+1-n} \frac{1-q^n}{1-aq}\,p_{n-1}(q^s;aq;q),
$$
$$
[(1-q^s)\,e^{-\partial_s}-aI]\,p_n(q^s;a;q)=(1-a)\,p_{n+1}(q^s;a/q;q),
$$
respectively.

To conclude this example let us mention that the Charlier
$q$-polynomials $c_n^{(\mu)}(s;q)$ on the lattice $x(s)=q^s$,
introduced in \cite{alar1}, are a  special case of the Wall polynomials, considered above.

\subsubsection{Discrete Laguerre $q$-polynomials
$L_n^{(\alpha)} (x;q)$}
Let us consider now the Hamiltonian, associated with the discrete Laguerre
$q$-polynomials. For these polynomials  $\sigma(s)=q^{-1}q^s$ and
$\sigma(s)+\tau(s)\nabla x_1(s)=aq^s(q^s+1)$ \cite{kost}.
If we put $A(s)=\sqrt{\nabla x_1(s)}$ and use the condition
\refe{first_cond_lin}, then  we obtain that $\alpha=1$ and
$\varsigma = \sqrt{\, q}$. Next we use the second condition
\refe{second_cond_lin}, but the corresponding expression is
a constant if and only if $a=q^{-1/2}$, i.e., not for
any value of the parameter $a$ the $q$-Hamiltonian, associated
with the discrete Laguerre $q$-polynomials, can be factorized
by using the $\alpha$-operators, which satisfy the corresponding
commutation relation. But as in the previous case,  one can define the functions
$$
\Phi_n^{(\alpha)}(s;q)=d_n^{-1}(\alpha) \frac{q^{s/2(\alpha+1)}}{\sqrt{(-q^s;q)_\infty}}\,
L_n^{(\alpha)} (x;q),\quad x:=x(s)=q^s,
$$
where the normalization constant $d_n(\alpha)$ is given by
$$
d_n(\alpha):=\frac{(-q^{\alpha+1},-q^{\alpha};q)_\infty^{1/2}}{(q;q)_\infty}
\frac{(q;q)_{n+\alpha}^{1/2}}{q^{n/2}(q;q)_n^{1/2}},
$$
and the discrete Laguerre $q$-polynomials $L_n^{(\alpha)} (x;q)$ are
$$
L_n^{(\alpha)} (x;q):=\frac1{(q;q)_n}\, {_2\phi_1}\left( \begin{array}{c} q^{-n}, x \\
0 \end{array} \Bigg\vert q;q^{n+\alpha+1}\right).
$$
These functions satisfy the discrete orthogonality relation  $\sum_{k=0}^\infty
\Phi_n^{(\alpha)} (k;q)\Phi_m^{(\alpha)} (k;q)=\delta_{n,m}$ and
$$
\g H^{(\alpha)} (s;q)\Phi_n^{(\alpha)} (s;q)=
\frac{1-q^n}{1-q}\Phi_n^{(\alpha)} (s;q),
$$
where the $q$-Hamiltonian is a difference operator of the form
$$
\g H^{(\alpha)} (s;q)=\frac{1}{1-q}\left\{\left[1+(1+q^{-\alpha})q^{-s}\right]I-q^{-\frac{s+\alpha}2} \left[\sqrt{1+q^s} \, e^{\partial_s}+
e^{-\partial_s} \sqrt{1+q^s}\right]\,q^{-\frac s2}\right\}.
$$

In this case the mutually adjoint operators,  that factorize the above $q$-Hamiltonian, are
$$
\begin{array}{l}
\dst \g a(s;\alpha;q) :=\frac1{\sqrt{1-q}}
\left(q^{-\frac\alpha2}I-e^{-\partial_s}\sqrt{1+q^s} \right)q^{-\frac s2}, \\[0.6cm]
\dst \g a^\dag(s;\alpha;q):= \frac1{\sqrt{1-q}}\,q^{-\frac s2}
\left(q^{-\frac\alpha2}I-\sqrt{1+q^s}\,e^{\partial_s} \right).
\end{array}
$$
They satisfy the commutation relation
$$
  \g a(s;\alpha-1;q) \g a^\dag(s;\alpha-1;q)  - q \g a^\dag (s;\alpha;q) \g a(s;\alpha;q)=I,
$$
from which it follows that
\[\begin{split}
&\g a(s;\alpha;q)\, \Phi_n^{(\alpha)} (s;q)=\sqrt{\frac{1-q^{n}}{1-q}}\, \Phi_{n-1}^{(\alpha+1)}(s;q),\\
&\g a^\dag (s;\alpha-1;q)\, \Phi_n^{(\alpha)} (s;q)=\sqrt{\frac{1-q^{n+1}}{1-q}}\, \Phi_{n+1}^{(\alpha-1)}(s;q).
\end{split}
\]
As is the previous case, these relations yield precisely the explicit form of the forward and backward shift operators,
respectively, for the $q$-Laguerre polynomials $L_n^{(\alpha)} (q^s;q)$ (see \cite{kost}, p. 109, (3.21.7) and (3.21,9)).


\subsubsection{Other cases in the $q$-Hahn tableau}

If one now applies the theorems \ref{el-teo-q-lam} and
\ref{el-teo-q-lin} to the families of big $q$-Jacobi polynomials with $b=0$
(big $q$-Laguerre), little $q$-Jacobi, $q$-Meixner, $q$-Kravchuck,
quantum $q$-Kravchuk, affine $q$-Kravchuk and  alternative
$q$-Charlier polynomials, then  for all these cases it is impossible to solve the problem 1.
The  big and little $q$-Jacobi polynomials do not admit a dynamical algebra because the corresponding
eigenvalues are not $q$-linear functions of $n$ and in the
other cases one of the two conditions of theorem
\ref{el-teo-q-lin} fails.

As an example let us consider the case of the $q$-Meixner polynomials.
For the $q$-Meixner polynomials \cite{kost}, we have $x(s)=q^{-s}$,
$\sigma(s)=c(x(s)-bq)/q$, $\sigma(s)+\tau(s)\nabla x_1(s)=(x(s)-1)(x(s)+bc)$
and $\lambda_n=q^{1/2}\frac{1-q^n}{(1-q)^2}$. Then, the condition
\refe{first_cond_lin} gives
$$
{q^{{\frac{1}{2}} - \alpha }}\sqrt{{{\frac{{q^{-2\alpha+1}}
\left({q^{\alpha}}-{q^s}\right)\left(bc{q^s}+{q^{\alpha }} \right)
\left( {q^{\alpha }}- b{q^{s+1}} \right) }{\left({q^s}-q\right)
\left(bc{q^s}+q \right)\left(b{q^{s+1}}-1 \right) }}}}=\q.
$$
After a careful study of the left-hand side one arrives at the
conclusion that it is constant if and only if
$\alpha=2/3$, $b=q^{-4/3}$, $c=-q^{5/3}$ or
$\alpha=2/3$, $b=q^{-5/3}$, $c=-q^{4/3}$. After substituting this in the second
condition \refe{second_cond_lin}, it becomes clear that the only possibility is the
first one, but it corresponds to a non-positive case.

Let us consider this case in more details.

The $q$-Meixner polynomials  (\cite{kost}, p.95),
\bq
M_n(q^{-x};b,c;q) := \, _2\phi_1\left(\begin{array}{c|c} q^{-n},q^{-x} \\
[-0.3cm] & q;\dst \frac{-q^{n+1}}{c} \\[-0.2cm]bq \end{array} \right)\,, \quad
0<b<q^{-1},\quad c>0\,, \label{(37)}
\eq
depend on two parameters $b$ and $c$ (in addition to the base $q$)
and occupy one level higher in the Askey scheme than $q$-Charlier
and Al-Salam--Carlitz polynomials of type II. They satisfy a difference
equation
$$
\left[B(x)\,(1-e^{\partial_x}) + D(x)\,(1-e^{-\partial_x})
\right]\,M_n(q^{-x};b,c;q)= (1-q^n)\,M_n(q^{-x};b,c;q)\,,$$
$$B(x) = c\,q^x\,(1-b\,q^{x+1}),  \quad D(x)= (1-q^x)(1+bc\,q^x).
$$
So $q$-Meixner functions, defined as
$$
\Phi_n^M(x;b,c;q):= d^{-1}_n(b,c)\,\left[\frac{c^x\,(bq;q)_x}{(q;q)_x\,
(-bcq;q)_x}\right]^{1/2}\,q^{x(x-1)/4}\,M_n(q^{-x};b,c;q)\,,
$$
are eigenfunctions of a difference ``Hamiltonian" $\g H^M(x;b,c;q)$,
$$
\g H^M(x;b,c;q)\Phi_n^M(x;b,c;q) = \frac{1-q^n}{1-q}\Phi_n^M(x;b,c;q)\,,
$$
where
$$
\g H^M(x;b,c;q):= \frac{1}{1-q}\,\left[B(x) + D(x) -
B^{1/2}(x)\, e^{\partial_x}\,D^{1/2}(x) -
D^{1/2}(x)\,e^{-\partial_x}\,B^{1/2}(x)\right]\,.
$$
The $q$-Meixner functions satisfy the discrete orthogonality
relation
$$
\sum_{k=0}^{\infty}\,\Phi^M_m(k;b,c;q)\,\Phi^M_n(k;b,c;q) = \delta_{mn}\,.
$$
One can factorize the ``Hamiltonian'' $\g H^M(x;b,c;q)$,
$$
\g H^M(x;b,c;q) = a^{\uparrow}_M(x;b,c;q)\,a^{\downarrow}_M(x;b,c;q),
$$
by means of the ``lowering'' and ``raising'' difference operators
\bq\begin{split}
\g a_M^{\downarrow}(x;b,c;q):= & \frac{1}{\sqrt{1-q}}\,\left[\,e^{\partial_x}\,D^{1/2}(x)-B^{1/2}(x)\,\right], \\
\g a_M^{\uparrow}(x;b,c;q):= & \frac{1}{\sqrt{1-q}}\,\left[\,D^{1/2}(x)\,e^{-\partial_x}-B^{1/2}(x)\,\right]\,.    \label{(44)}
\end{split}
\eq
The difference operators \refe{(44)} satisfy a $q$-commutation relation
of the form
$$
\g a_M^{\downarrow}(x;b/q,cq;q)\,\g a_M^{\uparrow}(x;b/q,cq;q)-
q\,\g a_M^{\uparrow}(x;b,c;q)\, \g a_M^{\downarrow}(x;b,c;q)= I\,.
$$
Their action on the $q$-Meixner functions is given by
\bq\ba{l}
\g a_M^{\downarrow}(x;b,c;q)\,\Phi^M_n(x;b,c;q)=
\sqrt{\frac{1-q^n}{1-q}}\,\Phi^M_{n-1}(x;bq,c/q;q),\\[4mm]
\g a_M^{\uparrow}(x;b/q,cq;q)\,\Phi^M_n(x;b,c;q)=  \sqrt{\frac{1-q^{n+1}}{1-q}}
\,\Phi^M_{n+1}(x;b/q,cq;q)\,,
\ea
\label{(46)}
\eq
that is, they not only lower and raise, respectively, the index $n$,
but alter also the parameters $b$ and $c$. The formulae
\refe{(46)} are equivalent to the statement that the forward and backward shift
operators for the $q$-Meixner polynomials \refe{(37)} have the form
(see [7], p.95, (3.13.2) and (3.13.8))
$$
(1 - e^{\partial_x})\,M_n(q^{-x};b,c;q) = \frac{1-q^n}{c(1-bq)}\,
q^{-x}\,M_{n-1}(q^{-x};bq,c/q;q)\,,$$
$$
\left[cq^x(1-bq^x) - (1-q^x)(1+bcq^x)\,e^{-\partial_x}\right]\,
M_n(q^{-x};b,c;q) 
$$
$$= cq^x(1-b)\,M_{n+1}(q^{-x};b/q,cq;q)\,.
$$

It remains only to remind the reader that when the parameter $b$
in \refe{(37)} vanishes, the $q$-Meixner polynomials $M_n(q^{-x};0,c;q)$
coincide with the $q$-Charlier polynomials (\cite{kost}, p.112)
\bq
C_n(q^{-x};c;q) := \, _2\phi_1\left(\begin{array}{c|c} q^{-n},q^{-x} \\
[-0.3cm] & q;\dst \frac{-q^{n+1}}{c} \\[-0.2cm] 0 \end{array} \right)\,.
\label{(48)}
\eq
In this case $B(x)=c\,q^x$ and $D(x)=1-q^x$. The appropriate formulae
for the $q$-Charlier polynomials \refe{(48)} are therefore easy consequences
of the corresponding formulae for the $q$-Meixner polynomials
\refe{(37)} with the
vanishing value of the parameter $b$.

\subsection{Askey-Wilson  cases}
Using the linearity of the difference equation \refe{eqdif} we
can consider, with not loss of generality, the following
lattice $x(s)=\frac12(q^s+q^{-s})$, for which $\mu=0$. Then,
$$
\sigma(s)=C q^{-2s}\prod_{i=1}^4 (q^s-q^{s_i})=
q^{-2s}\prod_{i=1}^4 (q^s-{z_i}),\qquad
\sigma(-s-\mu)=C q^{2s}\prod_{i=1}^4 (q^{-s}-z_i).
$$

It is well known that the general case when the zeros of
$\tilde\sigma$, namely $z_1 z_2 z_3 z_4 \ne 0$, corresponds to the
Askey-Wilson polynomials \cite{kost}.
We will use the theorem \ref{el-teo-q-lin} to solve the
factorization problem for the whole $q$-Askey tableau \cite{kost}
and Nikiforov-Uvarov tableau \cite{niuv2,atrasu}.

The Askey-Wilson polynomials on the lattice
$x(s)=\frac12(q^s+q^{-s})=\cos \theta$ where $q^s=e^{i\theta}$,
defined by \cite{kost}
$$p_n(x(s);a,b,c,d|q)=\frac{(ab,ac,ad;q)_n}{a^n}
\,_4\varphi_3 \left(\!\!\begin{array}{c}
q^{-n},\,q^{n-1}abcd,\, aq^{-s},\, aq^{s}\\[0.3cm]
ab,\,ac,\,ad\end{array}\Bigg|q;q \!\!\right),
$$
where $a={z_1}$, $b={z_2}$, $c={z_3}$, $d={z_4}$ and
$\mu=0$. Their orthogonality relation is of the form
$$
\int_{-1}^1 \omega (x) p_n(x;a,b,c,d)p_m(x;a,b,c,d)
dx =\delta_{nm}d_n^2,
$$
where
$$
\omega(x)=\frac{h(x,1)h(x,-1)h(x,q^{\half})h(x,-q^{\half})}
{2\pi\sqrt{1-x^2}\, h(x,a)h(x,b)h(x,c)h(x,d)},  \qquad
h(x,\alpha)=\prod_{k=0}^\infty [1-2\alpha x q^k+\alpha^2 q^{2k}],
$$
and the norm is given by
$$
d_n^2 = \frac{(abcdq^{n-1},abcdq^{2n};q)_\infty}{(q^{n+1},abq^n,acq^n,adq^n,bcq^n,bdq^n,cdq^n;q)_\infty}.
$$

The Askey-Wilson functions can be defined by
$$
\Phi_n(s;q)={\frac{\sqrt{\omega(s)} A(s)}
{d_n}} p_n(x(s);a,b,c,d),\quad x(s)=\cos\theta,
\quad q^s=e^{i\theta}.
$$

Taking $A(s)=\sqrt{\nabla x_1(s)}$, we have the orthogonality
$\int_{-1}^1 \Phi_n(s;q)\Phi_m(s;q)/\nabla x_1(s)dx=\delta_{n,m}$,
and $\g H(s;q)\Phi_n(s;q)=\lambda_n\Phi_n(s;q)$, where
$\lambda_n=q(q^{-n}-1)(1-abcdq^{n-1})$, and $\g H(s;q)$
is given by \refe{gen_hamilt},
$\sigma(s)= C_\sigma q^{-2s}(q^s-a)(q^s-b)(q^s-c)(q^s-d)$
and $\mu=0$. Thus the Hamiltonian, associated with these Askey-Wilson functions, is
$$
\begin{array}{c}
\dst \!\!\g H(s;q)\!=\!\frac{-1}{k_q^2\sqrt{\sin \theta}}\!
\left(\!\frac{\sqrt{\sigma(s)\sigma(-s+1)}}{\sin (\theta +\racion{i}{2}\log q)\sqrt{\sin (\theta+i\log q)}} e^{-\partial_s}+
\frac{\sqrt{\sigma(s+1)\sigma(-s)}}{\sin (\theta \!-\!\racion{i}{2}\log q)\sqrt{\sin (\theta-i\log q)}} e^{\partial_s}\!
\!\right) \\[0.6cm]
\dst +\frac{1}{k_q^2 \sin \theta}\left(\frac{\sigma(s)}{\sin (\theta +\racion i 2 \log q)}+
\frac{\sigma(-s)}{\sin (\theta -\racion i 2 \log q)}\right)I.
\end{array}
$$

If we now use the Remark \ref{rem-red-no-lin}, then the
first condition of the Theorem \ref{el-teo-q-lin} holds
for $\alpha=\half$. In fact, a straightforward calculations shows that
$$
\Delta x_1(s+\gamma)=
\frac{k_q}{2}q^{-s-\gamma}(q^{s+\gamma}+1)(q^{s+\gamma}-1).
$$
Hence, the condition \refe{first_cond_lin} has the form
 ($e^{i\theta}=q^s$)
$$\ba{l}
\dst q^{-2\alpha+1}\left(\frac{q^{s-\half}-1}{q^{s-\alpha}-1}
\sqrt{\frac{(q^{s-1}-1)(q^s-1)}
{(q^{s-\alpha-\half}-1)(q^{s-\alpha+\half}-1)}}
\sqrt{\prod_{i=1}^4
\frac{(q^{s-\alpha}-z_i)(q^{-s+\alpha}-z_i)}
{(q^{s}-z_i)(q^{-s+1}-z_i)}}
\right)\times\\[5mm]
\qquad\qquad\qquad\qquad \dst\frac{q^{s-\half}+1}{q^{s-\alpha}+1}
\sqrt{\frac{(q^{s-1}+1)(q^s+1)}
{(q^{s-\alpha-\half}+1)(q^{s-\alpha+\half}+1)}}=\q.
\ea
$$
If we look at the expression in the second line,  we
find that it is a constant if and only if
$\alpha=\half$, and thus the first condition transforms into
\begin{equation} \label{cond1_ask_wils}
\sqrt{\frac{\prod_{i=1}^4
(q^{s-1/2}-z_i)(q^{-s+1/2}-z_i)}
{\prod_{i=1}^4 (q^{s}-z_i)(q^{-s+1}-z_i)}}=\q.
\end{equation}
Then, the simplest case for which the condition
\refe{first_cond_lin} holds is when $z_1z_2=q^\half$ and
$z_3z_4=q^\half$ (with the corresponding permutations of the roots).
In this case $\varsigma=1$. Since $\sigma(s)=\sigma(-s+\half)$,
the second condition gives
$$
\frac{1}{\nabla x_1(s)} \left(\frac{\sigma(s+\half)}{\Delta x(s)}
+\frac{\sigma(-s+\half)}{\nabla x(s)}-
\frac{\sigma(s)}{\nabla x(s)}-\frac{\sigma(-s)}{\Delta x(s)} \right)=0.
$$

Thus, we have $a_{1/2}^{\uparrow}(s;q) a_{1/2}^\downarrow(s)=
\g H(s;q)$ and $[a_{1/2}^{\downarrow}(s;q),a_{1/2}^\uparrow(s)]_q=0$,
where
$$
\begin{array}{l}
\dst a_{1/2}^{\downarrow}(s;q)=
e^{\half \partial_s}\sqrt{\frac{\sigma(s)}{-k^2_q \sin \theta \sin (\theta +
\racion i 2 \log q)}} - e^{-\half \partial_s}\sqrt{\frac{\sigma(-s)}{-k^2_q
\sin \theta \sin (\theta -\racion i 2 \log q)}},\\[0.6cm]
\dst a_{1/2}^{\uparrow}(s;q)=\sqrt{\frac{\sigma(s)}{-k^2_q \sin \theta \sin
(\theta +\racion i 2 \log q)}}e^{-\half \partial_s}
- \sqrt{\frac{\sigma(-s)}{-k^2_q \sin \theta \sin (\theta -\racion i 2 \log q)}}
e^{\half \partial_s}.
\end{array}
$$
Since the $\alpha$-operators are commuting, this case is not so
interesting in applications (e.g., for $q$-models of the harmonic oscillators).
A special case of the Askey-Wilson polynomials are the
continuous $q$-Jacobi polynomials corresponding to the roots
$a=q^{\alpha'/2+1/4}$, $b=q^{\alpha'/2+3/4}$,
$c=-q^{\beta'/2+1/4}$, $d=-q^{\beta'/2+3/4}$ \cite{kost}, then we can
solve the problem 1 only in the case when $\alpha'=\beta'=-1/2$.\\
The next case is when one of the roots $z_i$ vanishes. This is the
0-Askey-Wilson polynomials or the continuous dual $q$-Hahn
\cite{kost}. In this case the first condition \refe{first_cond_lin}
(equivalently \refe{cond1_ask_wils}) holds only when
$$
(z_1,z_2,z_3)=(t,\racion12 - t,\racion 14), \qquad t\in \RR.
$$
With the above choice of $z_i$ it is impossible to fulfill the
second condition \refe{second_cond_lin} hence it is impossible to
obtain a simple closed dynamical algebra.

\subsubsection{Continuous $q$-Laguerre polynomials}

Let now consider the case when the Askey-Wilson
polynomials have two parameters equal to zero. In order that
the condition \refe{first_cond_lin} take place, the other
two non-vanishing parameters should satisfy that their product
is equal $q^\half$. Under this condition
$\varsigma=q^{-\half}$. Then
the second condition yields $\Lambda=
\frac{4C_\sigma(\sqrt{q}-1)}{k^2_q }$. Then, the Askey-Wilson
Hamiltonian with two zero roots of $\sigma$ admits a factorization
with a non-trivial dynamical algebra.
An example of this family is the continuous $q$-Laguerre polynomials
$P_n^{(a)}(x|q)$ \cite{kost}, $x(s)=\cos\theta$, when
$a=-1/2$ for which the Hamiltonian has the form ($A(s)=1$)
$$
\begin{array}{c}
\dst \!\!\g H(s;q)\!=\!-\frac{C_\sigma}{k_q\sin \theta}\left(
\frac{\sqrt{(q^{s+1}-1)(q^{s+1}-q^\half)(q^{-s}-1)(q^{-s}-q^\half)}}
{\sin (\theta +\racion{i}{2}\log q)}
e^{-\partial_s}+\right. \\[0.6cm]\dst\left.
\hspace{6cm}+
\frac{\sqrt{(q^{s}-1)(q^{s}-q^\half)(q^{1-s}-1)(q^{1-s}-q^\half)}}
{\sin (\theta -\racion{i}{2}\log q)} e^{\partial_s}
\right)+ \\[0.6cm]
\dst -\frac{C_\sigma}{k_q \sin \theta}\left(\frac{(q^{s}-1)(q^{s}-q^\half)}
{\sin(\theta+\racion i2\log q)}+\frac{(q^{-s}-1)(q^{-s}-q^\half)}
{\sin (\theta -\racion i 2 \log q)}\right)I,
\end{array}
$$
and
$$
\begin{array}{l}
\dst a_{1/2}^{\downarrow}(s;q)=\frac{\sqrt{-C_\sigma}}{k_q\sin\theta}\left(
e^{\half \partial_s}\sqrt{(q^{s}-1)(q^{s}-q^\half)}-
 e^{-\half \partial_s}\sqrt{(q^{-s}-1)(q^{-s}-q^\half)}\right),
\\[0.6cm]
\dst a_{1/2}^{\uparrow}(s;q)=
\frac{\sqrt{-C_\sigma}}{k_q\sin\theta}\left(
\sqrt{(q^{s}-1)(q^{s}-q^\half)}e^{-\half \partial_s}-
\sqrt{(q^{-s}-1)(q^{-s}-q^\half)} e^{\half \partial_s}\right).
\end{array}
$$
Then for the functions
$$
\Phi_n(x(s))=\sqrt{\frac{(q^\half;q)_n(q,q^\half;q)_\infty\omega(s)}
{(q;q)_n}}\,_3\varphi_2 \left(\!\!\begin{array}{c}
q^{-n},q^{-s}, q^{s}\\[0.3cm]
q^\half,0\end{array}\Bigg|q;q \!\!\right), \quad
q^s=e^{i\theta},\quad x(s)=\cos\theta,
$$
we have the orthogonality $\int_{-1}^1 \Phi_n(x;q)\Phi_m(x;q)dx=
\delta_{n,m}$, and
$$
\g H(s;q)\Phi_n(s;q)=q(q^{-n}-1)\Phi_n(s;q),\quad
\g H(s;q)=a_{1/2}^{\uparrow}(s;q) a_{1/2}^\downarrow(s),\quad
$$
and
$$[a_{1/2}^{\downarrow}(s;q),a_{1/2}^\uparrow(s)]_{q^{-1/2}}=
\frac{4C_\sigma(\sqrt{q}-1)}{k^2_q },
$$
Thus choosing $C_\sigma=-\frac{k^2_q}{4(1-\sqrt{q})}$ we
obtain the relation
$[a_{1/2}^{\downarrow}(s;q),a_{1/2}^\uparrow(s)]_{q^{-1/2}}=I$.

The case of Askey-Wilson polynomials with three zero parameters
is analogous to the case of one zero parameter and there is not
possible to solve the problem 1. An example of this case are the
continuous big $q$-Hermite polynomials \cite{kost}.

\subsubsection{Continuous $q$-Hermite polynomials}
Finally, if one takes the Askey-Wilson polynomials with vanishing parameters
$a,b,c,d$, this gives the continuous $q$-Hermite
polynomials \cite{kost}. In this case $\sigma(z)=C_\sigma q^{2z}$.

Let choose $A(s)=\sqrt{\nabla x_1(s)}$. Taking into account
that this family is a special case of the Askey-Wilson polynomials
when all parameters $a=b=c=d=0$, one directly obtains, by using
$\varsigma=1/q$, that in this case the
Hamiltonian is given by
$$
\begin{array}{c}
\dst \!\!\g H(s;q)\!=
\!\frac{-1}{k_q^2\sqrt{\sin \theta}}\! \left(
\!\frac{C_\sigma q}{\sin (\theta +\racion{i}{2}\log q)
\sqrt{\sin (\theta+i\log q)}} e^{-\partial_s}+
\frac{C_\sigma q}{\sin (\theta \!-\!\racion{i}{2}\log q)
\sqrt{\sin (\theta-i\log q)}} e^{\partial_s}\!\right)\! \\[0.6cm]
\dst+ \frac{1}{k_q^2 \sin \theta}\left(\frac{C_\sigma q^{2s}}
{\sin (\theta +\racion i 2 \log q)}+
\frac{C_\sigma q^{-2s}}{\sin (\theta -\racion i 2 \log q)}\right)I,
\end{array}
$$
and the $\alpha$-operators
$$
\begin{array}{l}
\dst a_{1/2}^{\downarrow}(s;q)=
e^{\half \partial_s}\sqrt{\frac{C_\sigma q^{2s}}{-k^2_q \sin \theta
\sin (\theta +\racion i 2 \log q)}}
- e^{-\half \partial_s}\sqrt{\frac{C_\sigma q^{-2s}}{-k^2_q
\sin \theta \sin (\theta -\racion i 2 \log q)}}\\[0.6cm]
\dst a_{1/2}^{\uparrow}(s;q)=
\sqrt{\frac{C_\sigma q^{2s}}{-k^2_q \sin \theta
\sin (\theta +\racion i 2 \log q)}}e^{-\half \partial_s}
- \sqrt{\frac{C_\sigma q^{-2s}}{-k^2_q \sin \theta \sin (\theta -\racion i 2 \log q)}}e^{\half \partial_s}
\end{array},
$$
are such that
$$
a^{\uparrow}(s;q)a^{\downarrow}(s;q)=\g H(s;q)\quad \mbox{ and }
\quad [a^{\downarrow}(s;q),a^{\uparrow}(s;q)]_{1/q}=
\frac{4C_\sigma}{k_q}.
$$
Notice that for getting the normalized commutation relations it
is sufficient to choose $C_\sigma=k_q/4$.

Another possible choice is $A(s)=1$ \cite{atsu2},
hence a straightforward calculus shows that the two
conditions in Theorem \ref{el-teo-q-lin} are true if
$\varsigma = q^{-1}$, thus $\Lambda={4C_\sigma}{k_q^{-1}}$.
With this choice the orthogonality of the functions $\Phi_n$
is $\int_{-1}^1 \Phi_n(s;q)\Phi_m(s;q)dx=\delta_{n,m}$.
In this case, the Hamiltonian is equal to
$$
\g H(s;q)=
\frac{C_\sigma q}{k_q^2}\left\{\frac{e^{-\partial_s}}{\sin
\theta \sin (\theta +\racion i 2 \ln q)}+
\frac{e^{\partial_s}}{\sin (\theta -\racion i 2 \ln q)\sin \theta}
- \frac{4}{\sqrt{q}}\left(1-\frac{1+q}{q+q^{-1}-2\cos
2\theta}\right) I \right\},
$$
and
$$
\begin{array}{l}
\dst \g a^{\downarrow}(s;q) := \g
a^{\downarrow}_{1/2}(s)=\frac{\sqrt{-C_\sigma}}{k_q \sin
\theta}\left(\ e^{\half \partial_s}q^s -
 e^{-\half \partial_s}q^{-s}\right), \\[0.6cm]
\dst \g a^{\uparrow}(s;q) := \g
a^{\uparrow}_{1/2}(s)=\frac{\sqrt{-C_\sigma}}{k_q \sin
\theta}\left( q^s  e^{-\half \partial_s}  - q^{-s}  e^{\half
\partial_s}  \right).
\end{array}
$$
In terms of these operators
$$
\g H(s;q)=a^{\uparrow}(s;q) a^{\downarrow}(s;q) \quad \mbox{ and } \quad
[a^{\downarrow}(s;q),a^{\uparrow}(s;q)]_{q^{-1}}=\frac{4C_\sigma}{k_q}.
$$
As before we now can choose $C_\sigma=k_q/4$.
This case was first considered in \cite{atsu1}, see also \cite{alatat}.

\appendix

\section*{Appendix}
Here we present a simple proof of the Theorem \ref{el-teo-q-lam}
when $\alpha=0$. Notice that in this case $b(s;q)=a^\dag(s;q)$, i.e.,
the operator $b$ is the adjoint of $a$.\\
\noindent\textbf{Theorem \ref{el-teo-q-lam}a}: Let
$\left\{\Phi_n(s;q)\right\}$ be the eigenfunctions of the operator
$\g H(s;q)$, corresponding to the eigenvalues
$\left\{\lambda_n(q)\right\}$. Suppose that the $\g H(s;q)$ admits
the factorization \refe{fac}. If the operators $a(s;q)$ and $b(s;q)$
in \refe{fac} satisfy the $q$-commutation relation
$[a(s;q)\,,\,b(s;q)]_q = I$, then the eigenvalues $\lambda_n(q)$ are
$q$-linear or $q^{-1}$-linear functions of $n$, i.e., either
$\lambda_n(q) = C_1\,q^n + C_3$ or $\lambda_n(q)= C_2\,q^{-n}+C_3$,
respectively.
\\
 \bdm Obviously, the operator $\g H(s;q)$ is diagonal
in the basis consisting of its eigenfunctions $\Phi_n(s;q)$. By
hypothesis, this operator admits at the same time the factorization
\refe{fac} in terms of $a(s;q)$ and $b(s;q)$. But according to
proposition \ref{pro1} the function $a(s;q)\Phi_n(s;q)$ is also the
eigenfunction of $\g H(s;q)$, associated with the eigenvalue
$q^{-1}[\lambda_n(q)-1]$. Hence the $a(s;q)$ is either the lowering
operator or the raising operator. In the former case this means that
$$
q^{-1}\,[\lambda_n(q)- 1] = \lambda_{n-1}(q)+ C \, ,
$$
from which it follows that $C_2 = 0$ (i.e., the spectrum
$\left\{\lambda_n(q)\right\}$ is a $q$-linear one) and
$C= q^{-1}\,[(1-q)\,C_3 - 1]$. In latter case the
corresponding relation is
$$
q^{-1}\,[\lambda_n(q)- 1] = \lambda_{n+1}(q) + C\, ,
$$
which holds when $C_1 = 0$ (i.e., the spectrum is $q^{-1}$-linear)
and $C = q^{-1}\,[(1-q)\,C_3 - 1]$. The proof of the theorem
is thus complete.
\edm

\noindent {\bf Acknowledgments:}
The research of RAN has been partially supported by the
Ministerio de Ciencias y Tecnolog\'{\i}a of Spain under
the grant BFM 2003-06335-C03-01, the Junta de Andaluc\'{\i}a
under grant FQM-262. The participation of NMA in this work has
been supported in part by the UNAM--DGAPA project IN102603-3
``\'Optica Matem\'atica'' and the Junta de Andaluc\'{\i}a
under grant 2003/2. Finally, we are grateful to both the referees for
their remarks and suggestions, which helped us to improve the
exposition of our results.

\bibliographystyle{plain}

\end{document}